\documentclass[a4paper,11pt]{amsart}
\usepackage{amssymb}
\usepackage[dvips]{graphicx}
\usepackage{amscd}
\newcommand{\comm}[1]{}

\def\dist{\operatorname{dist}}

\def\ti{\tilde}

\def\({\left(}
\def\){\right)}

\def\uli{\underline}

\def\raw{\rightarrow}

\def\no={\neq}
\def\sm{\setminus}

\def\C{{\mathbb C}}

\def\N{{\mathbb N}}

\def\P{{\mathbb P}}
\def\Q{{\mathbb Q}}

\def\BB{{\mathcal B}}

\def\MM{{\mathcal M}}
\def\NN{{\mathcal N}}
\def\OO{{\mathcal O}}

\def\al{\alpha}
\def\be{\beta}
\def\ga{\gamma}
\def\de{\delta}

\def\vep{\varepsilon}

\def\th{\theta}

\def\ka{\kappa}
\def\la{\lambda}

\def\om{\omega}

\def\La{\Lambda}

\theoremstyle{plain}
\newtheorem{Main}{Theorem}

\newtheorem{Thm}{Theorem}[section]
\newtheorem{Prop}[Thm]{Proposition}
\newtheorem{Lem}[Thm]{Lemma}

\newtheorem{Cor}[Thm]{Corollary}

\theoremstyle{remark}
\newtheorem{Rem}[Thm]{Remark}
\newtheorem{Def}[Thm]{Definition}

\begin{document}
\begin{center}
\end{center}
\title[Rational Misiurewicz maps are rare]
{Rational Misiurewicz maps are rare}
\author{Magnus Aspenberg}
\address{Universit\'e Paris-Sud, Laboratoire de Math\'ematique, B\^atiment 425, UMR 8628, 91405 Orsay, France}
\email{Magnus.Aspenberg@math.u-psud.fr}
\thanks{The author gratefully acknowledges funding from the Swedish Research Council}
\begin{abstract}
We show that the set of Misiurewicz maps has Lebesgue measure zero in the space of rational functions for any fixed degree $d \geq 2$.
\end{abstract}

\maketitle


\section*{Introduction}

Let $f(z)$ be a rational function of a given degree $d \geq 2$ on the Riemann sphere $\hat{\C}$. Let $Crit(f)$ be the set of critical points of $f$. Define
\[
P^k(f,c) =  \overline{ \bigcup_{n > k} f^n(c)} \text{ and } P^k(f) =  \bigcup_{c \in Crit(f)} P^k(f,c)  .
\]
Set $P^0(f,c)=P(f,c)$. The set $P(f)=P^0(f)$ is the {\em postcritical set} of
$f$. We will also use the
notion postcritical set for $P^k(f)$ for some suitable $k \geq 0$. Denote
by $J(f)$ the Julia set of $f$ and $F(f)$ the Fatou set of $f$. Recall that the $\om$-limit set $\om(x)$ of a point $x$ is the set of all limit points of $\cup_{n \geq 0} f^n(x)$. A periodic point $x$ with period $p$ is a {\em sink} if there is a neighborhood around $x$ which is mapped strictly inside itself homeomorphically under $f^p$. Hence sinks are attracting cycles which are not super-attracting. We proceed with the following definition.

\begin{Def}
A {\em Misiurewicz} map $f$ is a non-hyperbolic rational map that has no parabolic cycles or sinks and such that $\om(c) \cap Crit(f) = \emptyset$ for every $c \in Crit(f)$ not belonging to a super-attracting cycle.
\end{Def}

We prove the following.
\begin{Main}
The set of Misiurewicz maps has Lebesgue measure zero in the space of rational functions for any fixed degree $d \geq 2$. \label{miserere}
\end{Main}

The notion of Misiurewicz maps has its origin in the paper \cite{MM} by M. Misiurewicz.
In the quadratic family $f_a(x)=1-ax^2$, where $a \in (0,2)$, a Misiurewicz map is a non-hyperbolic map where the critical point $0$ is non-recurrent. D. Sands showed in \cite{DS} that the set of parameters $a \in (0,2)$ for which $f_a$ is Misiurewicz, has Lebesgue measure zero. In 2003, S. Zakeri showed in \cite{SZ} that the Hausdorff dimension of the set of Misiurewicz maps in the quadratic family $f_a$ is full, i.e. equal to $1$. Conjecturally, a similar statement holds in higher dimensions too.

There has been some variations on the definition of complex Misiurewicz maps. The maps studied in Misiurewicz original paper \cite{MM} are, among other things, (real) maps $f$ for which every critical point $c$ has that $\om(c) \cap Crit(f) = \emptyset$ and there are no sinks. In the complex case a Misiurewicz map has often been referred to as a map for which the critical points land on repelling periodic points (without possibly the critical points laying in super-attracting cycles, as for example the point at infinity for polynomials). In \cite{SvS}, S. van Strien studies Misiurewicz maps with a definition similar to the definition in \cite{MM}. In \cite{GKS}, a Misiurewicz map is a meromorphic map $f$ for which every critical point $c \in J(f)$ has the property that $\om(c) \cap Crit(f)=\emptyset$. It should not be hard to generalize Theorem \ref{miserere} to the more general definition in \cite{GKS}. However, the main construction will be the same as in this paper.

Let $SupCrit(f)$ be the set of critical points belonging to super-attracting cycles. For $B(z,r)= \{ w :  |w-z| < r \}$, let
\[
U_{\de} = \bigcup_{c \in Crit(f) \sm SupCrit(f)} B(c,\de).
\]
Let $\MM^d$ be the set of Misiurewicz maps of degree $d$. Define
\begin{align}
M_{\de,k} = \{ f \in \MM^d: P^k(f) \cap U_{\de} = \emptyset \}
\label{dek-misse}
\text{ and } M_{\de} = \cup_{k \geq 0} M_{\de,k}.
\end{align}
If $f \in M_{\de,k}$ then we say that $f$ is {\em $(\de,k)$-Misiurewicz}. If $f \in M_{\de}$ then we say that $f$ is {\em $\de$-Misiurewicz }.
So every Misiurewicz map $f$ has some $\de > 0$ and $k \geq 0$ for which $P^k(f) \cap U_{\de} = \emptyset$.

Every rational map of degree $d$ can be written in the form
\begin{equation}
R(z)=\frac{P(z)}{Q(z)}=\frac{a_0+a_1z+\ldots+a_dz^d}{b_0+b_1z+\ldots+b_dz^d},
\label{rat}
\end{equation}
where $a_d$ and $b_d$ are not both zero. Without loss of generality we
may assume that $b_d = 1$. The case $a_d \neq 0, b_d = 0$ is treated
analogously. Hence, the set of rational functions of degree $d$ is a
$2d+1$-dimensional complex manifold and subset of the projective space
$\P^{2d+1}(\C)$. Now, simply take the measure on the coefficient space
in one of the two charts $a_d=1$ or $b_d=1$. There also is a
coordinate independent measure on the space of rational maps of a
given degree $d$, induced by the Fubini-Study
metric (see \cite{FS}). The Lebesgue measure on any of the two charts
is mutually absolutely continuous to the Fubini-Study measure.

A family of rational maps $R_a$ for $a \in V \subset \C^m$, where $V$ is open
and connected, is {\em normalized} if any two functions $R_a$ and $R_b$, $a,b
\in V$, are conformally conjugate if and only if $a = b$. If $f$ and $g$ are conformally conjugate then they are conjugate by a M\"obius transformation
\[
M(z) = \frac{\al + \be z}{\ga + \de z}.
\]
The set of M\"obius transformations forms a $3$-dimensional complex manifold. Introduce an equivalence relation $\sim$ on the parameter space, saying that $f \sim g$ if and only if $f = M^{-1} \circ g \circ M$, for some M\"obius transformation $M$. Every equivalence class is a complex $3$-dimensional manifold. These manifolds form a foliation of the space of rational functions of degree $d$ (see e.g. Frobenius Integrability Theorem in \cite{MiSp}). Hence to prove Theorem \ref{miserere}, by Fubini's Theorem, it suffices to consider families of normalised maps. If fact, we will consider $1$-dimensional slices of such maps. The following theorem is the main object of this paper.
\begin{Main} \label{discthm}
Assume that $R_a$, $a \in \C$, is an analytic normalized family of rational maps in a
neighborhood of $a=0$ and that $R_0$ is a Misiurewicz map.
Then the Lebesgue density at $a=0$ is strictly less than $1$ in the
set of $(\de,k)$-Misiurewicz maps for any $\de > 0$ and $k \geq 0$.
\end{Main}

Theorem \ref{miserere} then follows from Theorem \ref{discthm} and Fubini's Theorem.
Indeed, by Fubini's Theorem, the Lebesgue density of
the set of $(\de,k)$-Misiurewicz maps in the full parameter
space $\C^{2d+1}$, is strictly less than $1$ at the point
$R=R_0$. Since this is true for every Misiurewicz map $R$ and every
$\de > 0$ and $k \geq 0$, and since a set of positive Lebesgue measure must have Lebesgue density equal to $1$ almost everywhere, it follows that the
Lebesgue measure of the set of Misiurewicz maps of a given degree is
\[
m\biggl( \bigcup_{n,k \in \N} M_{1/n,k} \biggr) \leq \sum_{n,k \in \N} m(M_{1/n,k}) = 0.
\]

There are some similarities between the methods in this paper and the paper \cite{DS} by D. Sands. The existence of a continuation
of the postcritical set in the real case in \cite{DS} is replaced by a
similar idea, namely that of a holomorphically moving postcritical set
in the complex case. However, to prove the existence of such a set we will use results by S. van Strien in \cite{SvS}. This paper uses much of the ideas in \cite{MA} and some fundamental results from the paper by Benedicks-Carleson in \cite{BC2} (and \cite{BC1}).

\subsection*{Acknowledgements}
I am grateful to Michael Benedicks for many helpful discussions and
comments on, among other things, the transversality condition.
I am grateful to Jacek Graczyk, Duncan Sands, Nicolae Mihalache, Neil Dobbs for many useful remarks on a preliminary version. I am thankful for very interesting discussions with Dierk Schleicher. Finally, I am thankful to Nan-Kuo Ho for useful comments and encouragement.

This paper was written at the Department of Mathematics at Universit\'e Paris-Sud, whose hospitality I owe my thanks to.
\section{Some definitions}

Let $R_0(z)=R(z)=P(z)/Q(z)$ be the starting unperturbed rational map of degree $d=\max{(\deg(P),\deg(Q))}$
and assume that $R(z)$ is Misiurewicz. This means that for some $\de > 0$ we
have $U_{\de} \cap P^k(R) = \emptyset$ for some $k \geq 0$. In
addition, choose $\de$ such that $U_{2\de} \cap P^k(R) = \emptyset$.

The critical points $c_i$ may of course split into several critical
points under perturbations if the multiplicity of $c_i$ is higher than
$1$. For a detailed description of this phenomena we refer to \cite{MA},
Section 1.3 and the standard theory of zero sets of irreducible
polynomials \cite{BK}, Theorem 1, p. 386. Either some $c_i$ split into
several critical points $c_{ij}(a)$, which is analytic in $a$, where
$c_{ij}(a)$ is of lower multiplicity for $a \neq 0$ close to $0$, or
$c_j(a)$ moves holomorphically itself with constant multiplicity. The
latter case will be referred to as the non-degenerate case and the
other case as the degenerate case. In the degenerate case the critical
points $c_{ij}(a)$ emerging from $c_i$ form so called {\em critical
  stars}, so that $c_{ij}(a) \neq c_{ik}(a)$ if $j \neq k$, for $a
\neq 0$ and $a \in B(0,r)$.


We will study a critical point $c=c(a)$ dependent on the parameter
$a \in B(0,r)$, for some $r > 0$. We sometimes write $R(z,a) =
R_a(z)$. Put
\[
\xi_{n,j}(a)=R^n(v_j(a),a),
\]
where $v_j(a)=R(c_j(a),a)$ is a non-critical critical value and
$c_j(a) \in Crit(R_a)$.
(A priori there can be finite chains of critical
points mapped onto each other. Therefore we assume that $v_j(a)$ is
the last critical value).
For simpler notation, we sometimes drop the index $j$ and write only $\xi_{n,j}(a) = \xi_n(a)$.

We also make the following convention. Chosen $\de > 0$, we always
assume that the parameter disk $B(0,r)$ is chosen so that the critical
points $c_i(a)$ moves inside $B(c_i,\de^{10})$ as $a \in B(0,r)$.

Important constants are the real positive numbers $\de',\de''$ where $\de' \gg \de''$ shall only depend on the unperturbed function $R_0=R$ and slightly on the perturbation $r > 0$.
There are three constants $S,S_1$ and $\ti{N}$ defined implicitly in Lemmas \ref{initdist}, \ref{distortion2} and \ref{finret} respectively.

We use the spherical metric and the spherical derivative unless otherwise stated.

\begin{Def}
Given two complex numbers $A$ and $B$, we write $A \sim B$
meaning that there is a constant $C > 0$ depending only possibly on the unperturbed function $R$, $\de'$, and the perturbation
$r$ and such that the following holds:
\[
\biggl| \frac{A}{B} - 1 \biggr| \leq C.
\]
Moreover, we require that for any $\vep > 0$ there exist $\de',r > 0$ such that $C \leq \vep$.

If $A$ and $B$ are real and positive and $A \geq B - C$, then we write $A \gtrsim B$, if in addition for any $\vep > 0$ there exist $\de',r > 0$ such that $C \leq \vep A$. In particular $B \gtrsim A$ and $A \gtrsim B$ implies $A \sim B$.

\end{Def}

We will for simpler notation use the same $C$ for different
constants, even in the same series of equations. So for example
expressions like $C = 2C$ may appear.

\section{Expansion}

A fundament in getting distortion estimates is to have exponential growth of $|\xi_n'(a)|$.  The strategy we employ is similar to the one used in
\cite{MA}, Chapter 4, but the major difference here is that the postcritical
set is not (necessarily) finite. Recall that in \cite{MA} the critical points
for the unperturbed function all land on repelling periodic orbits. Moreover,
to prove that Misiurewicz maps have measure zero we need a certain
transversality condition (see below), which is not needed in \cite{MA}. The
transversality condition in the postcritically finite case means that the
critical value and the holomorphically moving repelling periodic orbit, on
which the critical point lands for the unperturbed function, do not coincide
for small perturbations. We will use Ma\~n\'e's Theorem and a result by S. van Strien to resolve this.

Recall that a compact set $\La$, which is invariant under $f$, is {\em
hyperbolic} if there are constants $C > 0$ and $\la > 1$ such that for any $z
\in \La$ and any $n \geq 1$,
\[
|(f^n)'(z)| \geq C \la^n.
\]
Equivalently, there is a metric $\varphi(z)$ which is expanding, meaning that
\begin{equation}
\varphi(f(z)) |f'(z)| > \varphi(z),
\label{Phimetrik}
\end{equation}
for all $z \in \La$.

The main result which we will use by Ma\~n\'e (see \cite{RM}) is the following.

\begin{Thm}[Ma\~n\'e's Theorem I] \label{mane}
Let $f:\hat{\C} \mapsto \hat{\C}$ be a rational map and $\La \subset J(f)$ a
compact invariant set not containing critical points or parabolic points. Then either $\La$ is a hyperbolic set or $\La \cap \omega(c) \neq \emptyset$
for some recurrent critical point $c$ of $f$.
\end{Thm}
\begin{Thm}[Ma\~n\'e's Theorem II] \label{mane2}
If $x \in J(f)$ is not a parabolic periodic point and does not intersect $\om(c)$ for some recurrent critical point $c$, then for every $\vep > 0$, there is a neighborhood $U$ of $x$ such that

\begin{itemize}
\item For all $n \geq 0$, every connected component of $f^{-n}(U)$ has diameter $\leq \vep$.

\item There exists $N > 0$ such that for all $n \geq 0$ and every connected component $V$ of $f^{-n}(U)$, the degree of $f^n |_V$ is $\leq N$.

\item For all $\vep_1 > 0$ there exists $n_0 > 0$, such that every connected component of $f^{-n}(U)$, with $n \geq n_0$, has diameter $\leq \vep_1$.

\end{itemize}
\end{Thm}
An alternative proof of Ma\~n\'e's Theorem can also be found by L. Tan and M. Shishikura in \cite{ST}.

A corollary of Ma\~n\'e's Theorem II is that a Misiurewicz map cannot have any Siegel disks, Herman rings or Cremer points (see \cite{RM} or \cite{ST}).
\comm{
Indeed, take a point $x$ on the boundary of a Siegel disk or Herman ring, call it $B$, and let $U$ be a neighborhood of $x$ as in Ma\~n\'e's Theorem II. Then there is some $y \in B \cap U$, with $\dist(y,\partial B) \geq r$, for some $r > 0$. Since a Misiurewicz map has no recurrent critical points, we get that every connected component $V$ of $f^{-n}(U)$ has diameter $d_n \raw 0$ as $n \raw \infty$. Now let $V$ be such a component intersecting the boundary of $B$ and containing $f^{-n}(y)$. Inside $B$, $f$ or an iterate of $f$ is conjugate to a rotation of angle $\th \notin \Q$. Hence there is a well defined inverse $f^{-1}$ (or $f^{-p}$) inside $B$. This implies that for every $\vep_1$ there is an $n_0 > 0$ such that we can find a $n \geq n_0$ such that $|f^{-n}(y)-y| \leq \vep_1$. But then since $V=f^{-n}(U)$ contains both $f^{-n}(x) \in \partial B$ and $f^{-n}(y) \in B$ we have that the diameter of $V$ must be greater than, say, $r/2$, if $\vep_1 < r/4$. The claim follows.
}
In particular, a Misiurewicz map has no indifferent cycles.

S. van Strien shows in \cite{SvS}:
\begin{Thm} [van Strien] \label{vanStrien}
Assume that $R_a$, for $a \in \C^m$ near the origin, is a
normalized family of rational maps and that $R_0$ is a Misiurewicz map, with non-degenerate critical points $c_1,\ldots,c_n$. Then there exists a neighborhood $W$ of $a=0$ and a holomorphic motion $h: P(R) \times W \raw \hat{\C}$ on the
postcritical set $P(R_0)$ of $R=R_0$ such that (writing $h_a(z)=h(z,a)$) for $a \in W$ the set $X_a = h_a(P(R_0))$ is mapped into itself by $R_a$ and
\[
h_a \circ R_0(z) = R_a \circ h_a(z), \text{ for all $z \in X_0$}
\]
Moreover, the map
\[
G(a) = (R_a(c_1(a))-h_a(R(c_1)), \ldots, R_a(c_n(a)) - h_a(R(c_n))
\]
is a local immersion for $a$ close to $0$.
\end{Thm}
Note that in our case we have $m=1$. We will use a generalization of
the uniqueness part of Thurstons Theorem in \cite{SvS} for Misiurewicz
maps. The holomorphically moving postcritical set in Theorem \ref{vanStrien}
plays the role of the (holomorphically moving) repelling periodic
orbits in \cite{MA}, on which the critical points land for the
unperturbed function.

For at least one $j$, we will need the following transversality condition,
namely that
\begin{equation} \label{transv}
\xi_{0,j}(a) \neq \mu_{0,j}(a), \text{ for $a \neq 0, a \in B(0,r)$},
\end{equation}
where $\mu_{n,j}(a) = h_a \circ R^n(v_j)$, $v_j=R(c_j)$.

In order to show Theorem \ref{vanStrien}, van Strien proves a generalized version of Thurston's Theorem, namely the following (see Theorem 3.3 in
\cite{SvS}).
\begin{Thm}
Assume that $R_1$ and $R_2$ are quasiconformally conjugate Misiurewicz maps and not Latt\'es maps. Then $R_1$ and $R_2$ are conjugate by a M\"obious
transformation.
\end{Thm}

The following lemma is heavily inspired by Theorem 3.3 in \cite{SvS}. 
\begin{Lem} \label{pullback}
If $\xi_{0,j}(a) = \mu_{0,j}(a)$ for all $j$, some $a \neq 0$, $a \in B(0,r)$, and such that any super-attracting cycle for $R_0$ persists under perturbations in $B(0,r)$, then $R_0$ and $R_a$ are quasi-conformally conjugate with bounded dilatation.
\end{Lem}
\begin{proof}
Note that the assumption implies that the critical points $c_j$ do not split under perturbation. Let $\La = P(R)$. Construct a holomorphic motion $h: B(0,r) \times \La \raw \hat{\C}$, (see for example \cite{MS}, Theorem III.1.6), such that it gives is a quasiconformal conjugacy $h_a$ on $\La$ (cf. Theorem 2.3.).

Now we want to extend this conjugacy to the whole sphere using a standard pullback argument. We extend $h_a$ by the $\la$-lemma (see \cite{MSS}) to the whole Riemann sphere. Call this extended motion $\ti{h}_a$. If the Julia set is not the whole Riemann sphere then there exist superattracting cycles, around which we can form a small neighborhood $N$ and get a conformal conjugacy in $N$ using the B\"ottcher coordinates. Then glue this conjugacy together with $\ti{h}_a$, so as to obtain a new $K_0$-quasiconformal homeomorphism $H_0$ which conjugates $R$ and $R_a$ on $\La \cup N$:
\begin{equation}
R_a \circ H_0(z) = H_0 \circ R_0(z), \text{ for $z \in \La \cup N$}. \label{conjugacy}
\end{equation}
If the Julia set is the whole Riemann sphere then take $H_0 = \ti{h}_a$ and $H_0$ is a conjugation on $\La$. In both cases write $Z$ as the set where the conjugation $H_0$ is valid and proceed as follows. Observe that the postcritical set $P(R)$ can be assumed to be contained in $Z$.
First we construct a sequence of homeomorphisms $H_n$, equal to
$H_{n-1}$ on $R_0^{-n}(Z)$, by pullback:
\[
H_n \circ R_0 = R_a \circ H_{n+1}.
\]
Since $H_n$ is a homeomorphism and $R_0$ and $R_a$ are covering maps on the
complement of the critical points which have the same corresponding multiplicities, the existence of the lifting follows by the General Lifting Lemma (see e.g. \cite{JaMu}, p. 390, Lemma 14.2).


It is also clear that the quasi-conformality of $H_n$ all have the same
upper bound $K$. Moreover, every map $H_0^{-1} \circ H_n$ fixes $Z$,
so the sequence $H_0^{-1} \circ H_n$ is equicontinuous. It follows that the family $H_n$ is equicontinuous. Therefore there is a subsequence which converges on compact subsets of $\hat{\C}$. But since also $H_n$ converges on $\cup_n R_0^{-n} Z$, to a conjugacy between $R_0$ and $R_a$, and since $\cup_n
R_0^{-n}(Z)$ is dense in $\hat{\C}$, $H_n$ must converge uniformly to a $K$-quasiconformal conjugacy on the whole sphere.

\end{proof}
If the function $G(a)$ in Theorem \ref{vanStrien} does not have an isolated
zero at $a=0$, then all critical points $c_j(a)$, for which $c_j(0) \in J(R_0)$, will still satisfy $\om(c_j(a_n)) \cap Crit(R_{a_n}) = \emptyset$, for a sequence $a_n \in B(0,r)$, where $a_n \raw 0$. If there is a super-attracting cycle, which bifurcates into a sink (non super-attracting attracting cycle), then $a=0$ is an isolated Misiurewicz point, since Misiurewicz maps have no sinks. If the super-attracting cycles persist, then by the above lemma there would be a sequence of maps $R_{a_n}$ which are quasiconformally conjugate Misiurewicz maps and $a_n \raw 0$ (cf. the beginning of the proof of Theorem 3.3 in \cite{SvS}). This would imply that the family $R_a$, $a \in B(0,r)$ is not normalized. Hence (\ref{transv}) must hold for at least one index $j$.

Now, if some critical point $c_i$ is degenerate, hence splits into
several critical points $c_{ij}(a)$ under perturbation, then we have
$c_{ij}(a) \neq c_{ik}(a)$ for $j \neq k$. So indeed, for at least one
index $j$ we must have $c_{ij}(a) \neq \mu_i(a)$.

For this reason, for any of the above two cases, by $\xi_n(a)$ we mean
$\xi_{n,j}(a)$ for this particular $j$ further on, unless otherwise
stated. Moreover, let us assume from now on that every super-attracting cycle for $R_0$ persists in the family $R_a$, $a \in B(0,r)$.

\subsection{Expansion near the postcritical set.}

By Ma\~n\'e's Theorem, the Misiurewicz
condition gives rise to expansion of the derivative in a (closed)
neighborhood of the postcritical set. More precisely, the
postcritical set $P^k(R)$ for a Misiurewicz map $R(z)$ is hyperbolic for some $k > 0$. Write
$P^k(R)=P$. Hence there exists a neighborhood $\NN$ of $P$ on which we have expansion in some metric
$\varphi(z)$. In other words,
\begin{equation}
\varphi(R(z)) |R'(z)| > \varphi(z),
\label{phi}
\end{equation}
for all $z \in \NN$. Thus, for some $C_0 > 0$ and $\la_0 > 1$,
\begin{equation}
|(R^j)'(z)| \geq C_0 \la_0^j,
\end{equation}
whenever $R^k(z) \in \NN$ for $k = 0,1,\dots,j$.

Let $\NN$ be so that $U_{10\de} \cap \NN = \emptyset$ and assume that $\NN$ is closed. Rewriting equation (\ref{phi}) with parameters, and shrinking $\NN$ if necessary, we get
\begin{equation}
\varphi(R(z,a)) |R'(z,a)| > \varphi(z),
\label{phi2}
\end{equation}
if $z \in \NN$ and $a \in B(0,r)$. By the compactness of $\NN$ and the continuity of equation (\ref{phi2}) with respect to $a$, $|(R^n)'(z,a)|$ grows exponentially for small $a$ whenever $z$ is inside $\NN$, i.e. we have the following lemma.
\begin{Lem} \label{expinside}
There exists some $\la > 1$, such that whenever $R_a^k(z) \in \NN$ for $k = 0,1,\dots,j$ and $a \in B(0,r)$, we have
\begin{equation}
|(R_a^j)'(z)| \geq C \la^j.
\end{equation}
\end{Lem}
Now choose the constant $\de' > 0$ so that
$
\{ z : \dist(z,P) \leq 10\de' \} \subset \NN.
$
Further, there will be more conditions on $\de'$ in Lemmas \ref{distortion1} and \ref{smallep} (so that we might have to diminish $\de'$). Define
\[
P_{\de'} = \{ z : \dist(z,P) < \de' \}.
\]

\section{Distortion Lemmas}
In this section we prove several distortion lemmas which will give control of $\xi_n(D_0)$ for dyadic disks $D_0 \subset B(0,r)$. We begin with a fundamental result, which says that the parameter and space derivatives are comparable as long as the space derivative grows exponentially. This idea was first introduced by Benedicks-Carleson in \cite{BC2}. The statement and proof here will be similar to the corresponding Proposition 4.3 in \cite{MA}. Let
\[
Q_n(a) = \frac{\partial R^n(v(a),a)}{\partial a} \bigg/
\frac{\partial R^n(v(a),a)}{\partial z},
\]
where $v(a)=R(c(a),a)$, and $v(a)$ is not a critical point. For the
unperturbed function ($a=0$) we write simply $Q_n(0)=Q_n$.

Write
\[
x(a) = \xi_0(a) - \mu_0(a),
\]
where $\mu_n(a) = h_a(R_0^n(v(0)))$, and $h_a$ is the holomorphic motion in Theorem \ref{vanStrien}. Observe that $x(a)$ is assumed to be not identically zero. This means that there is a $1 \leq k < \infty$ such that
\begin{equation} \label{xK1}
x(a) = K_1 a^k + \ldots.
\end{equation}
Given $\de' > 0$ and $r > 0$, in the following we let $\de'' > 0$ be a small fixed number, always satisfying $\de'' \ll \de'$ and $|\log \de''| \ll |\log |x(r)||$ (this holds if the  perturbation $r > 0$ is small enough).


In the following version of Proposition 4.3 in \cite{MA}, $\uli{\ga} = \log \la/2$, where $\la$ is as in Lemma \ref{expinside}.
\begin{Prop}
For every $\de > 0$ and sufficiently small $\de'' > 0$ there is an $r > 0$ such that the following holds. Assume that the parameter $a \neq 0$, $a \in B(0,r)$, satisfies $|\xi_n(a)-c(a)| \geq \de$ for all critical points $c(a)$ and all $n \geq 0$, and that
\begin{align}
|\partial R^n(v(a),a)/\partial z| &\geq Ce^{\ga n}, \quad \textrm{for}
\quad n=0,\ldots,m,  \\
|\partial R(\xi_n(a),a)/\partial z| &\leq B, \quad \forall n > 0,
\end{align}
where $\ga \geq \uli{\ga}$.
Then if $N > 0$ satisfies $|\xi_N(a)-\mu_N(a)| \geq \de''$, we have
for $n=N,\ldots m, $
\begin{equation}
|Q_n(a)-Q_N(a)| \leq \frac{1}{1000}|Q_N(a)|. \label{d}
\end{equation}
\label{da/dz}
\end{Prop}
\begin{Rem}
The number $N$ in the above proposition will be chosen for future use as the smallest possible, for which $|\xi_N(a)-\mu_N(a)| \geq \de''$ holds. This condition comes from one of the subsequent lemmas (Lemma \ref{xinprim}), where we have to "wind up" the comparison between $\xi_n'(a)$ and $(R_a^n)'(\mu_0(a))$ first before their ratio becomes stable. The Distortion Lemma then allows switching from $(R_a^n)'(\mu_0(a))$ to $(R_a^n)'(v(a))$.
\end{Rem}
\begin{Rem}
In \cite{MA} the proposition is stated in a stronger form. There we assume that the assumption $|\xi_n(a)-c(a)| \geq \de$ is replaced by $a \in \BB_{m-1,l}'$, which is a specific approach rate condition called the {\em basic assumption} (the idea was first introduced in \cite{BC1}, and \cite{BC2} ), meaning that
\[
|\xi_n(a) - c(a)| \geq e^{-\al n}
\]
for all critical points $c(a)$ and all $n \geq 1$. In our case this condition is much stronger than we need, since we only need to have (\ref{d}) fulfilled until the first return.
\end{Rem}

The use of the above lemma is similar to \cite{MA}.  Instead of having a repelling fixed point (or periodic orbit) we now have a hyperbolic set $\La = P$. In a neighborhood $\NN$ to this set the expansion of the derivative will be the fundament for the distortion estimates. Instead of looking at the whole parameter disk $B(0,r)$, we pick a smaller disk $D_0=D(a_0,r_0) \subset B(0,r)$ such that $|a_0|$ is significantly greater than $r_0$. This smaller disk $D_0$ will be mapped by $\xi_n$ onto a greater disk of a certain size $S = \OO(\de')$, still inside $\NN$. When this is achieved, we will get an important distortion estimate. There is a strong analogy with the ideas in \cite{MA}, where we use real analytic parameter families and we start with a small interval $[0,a_0]$, and pick a smaller interval $\om_0$ at the right end of $[0,a_0]$.

Before proving Proposition \ref{da/dz} we need some lemmas.
We have
\[
R_a(z) - R_a(w) = R_a'(w)(z-w) + \vep(w,z-w),
\]
where the error $\vep(w,z-w)$ depends on both the position of $w$ and $z-w$ (it also depends analytically on the parameter, but here skip that index for simpler notation). Since $\vep(w,z)$ is analytic in $z$ and $w$, and $\vep(w,z) = o(z)$ for fixed $w$, $z=0$ is a removable singularity, and we have
\[
\vep(w,z) = c_l z^l + \ldots
\]
for some $l \geq 2$. Recall that $|R_a^N(z)-R_a^N(w)| \geq \la |z-w|$ for some $\la > 1$ for all $z,w \in \NN$. For simplicity assume that $N=1$. We want to estimate the total error $E_1$ in
\[
R_a^n(z) - R_a^n(w) = (R_a^n)'(w)(z-w) + E_1.
\]
Let us state this as a lemma:
\begin{Lem} \label{smallep}
For every $\vep > 0$ there is a $\de' > 0$ and an $r > 0$ such that
the following holds. Fix some parameter $a \in B(0,r)$. If $R_a^j(z),
R_a^j(w) \in \NN$ and $|R_a^j(z)-R_a^j(w)| \leq \de'$ for
$j=0,\ldots,n$, $z \neq w$, then
\[
\biggl| \frac{R_a^n(z) - R_a^n(w)}{(R_a^n)'(w)(z-w)} - 1 \biggr| \leq \vep.
\]
\end{Lem}

\comm{
\begin{proof}[no 1]
By Lemma \ref{expinside} we have $|(R_a^n)'(z)| \geq C_0 \la_0^n$, so Lemma \ref{K} implies that the curvature of $R_a^n(\ga(t))$ is uniformly bounded by some constant, where $\ga(t) =z(1-t) + tw$, for $0 \leq t \leq 1$. Now, the length of the curve $R_a^n(\ga(t))$ is
\[
L = \int_0^1 |(R_a^n)'(\ga(t))||\ga'(t)|dt = |z-w| \int_0^1 |(R_a^n)'(\ga(t))|dt
\]

\end{proof}
(The lemma can also can be proven by induction over $n$. In this case we do not need the geometry estimates in the previous section. However, the proof becomes much longer).
}

\begin{proof}
The proof goes by induction over $n$. It is obviously true for $n=1$. Assume that it is true for $n-1$.

With $z_i=R_a^i(z), w_i=R_a^i(w)$ and $\vep(w_i,z_i-w_i)=\vep_i$ we have
\begin{align}
z_n-w_n &= R_a'(w_{n-1})(z_{n-1}-w_{n-1}) + \vep_{n-1} \nonumber \\
&= R_a'(w_{n-1})(R_a'(w_{n-2})(z_{n-2}-w_{n-2})+\vep_{n-2}) + \vep_{n-1} \nonumber \\
&= \ldots = (R_a^n)'(w_0)(z_0-w_0) + \sum_{j=1}^n \vep_{j-1} \prod_{i=j}^{n-1} R_a'(w_i) .\nonumber
\end{align}
Therefore,
\begin{equation}\label{above}
z_n-w_n = (R_a^n)'(w_0) \bigr( (z_0-w_0) + \sum_{j=1}^n \vep_{j-1} \prod_{i=0}^{j-1} \frac{1}{R_a'(w_i)} \bigl).
\end{equation}
We want to estimate the second term in (\ref{above}), which we call $E_2$. Using the induction assumption twice and the fact that $R_a$ is expanding on $\NN$ we obtain,
\begin{align}
|E_2| &\leq C \sum_{j=1}^n \frac{|z_{j-1}-w_{j-1}|^2}{|(R_a^j)'(w_0)|}
\leq  C \sum_{j=1}^n \frac{|(R_a^{j-1})'(w_0)|^2|z_{0}-w_{0}|^2}{|(R_a^j)'(w_0)|} \nonumber \\
&\leq C \sum_{j=1}^n \frac{|(R_a^{j-1})'(w_0)||z_{0}-w_{0}|^2}{|R_a'(w_{j-1})|} \nonumber \\
&\leq  C \sum_{j=1}^n |z_{j-1}-w_{j-1}||z_0-w_0| \nonumber \\
&\leq C |z_{n-1}-w_{n-1}||z_0-w_0| \leq C(\de') |z_0-w_0|.\nonumber
\end{align}
Thus if the maximum size $\de'$ of $|z_{n-1}-w_{n-1}|$ is bounded suitably, the lemma follows.


\end{proof}
\begin{Lem} \label{xinprim}
For every $\vep > 0$, if $\de' > 0$ is sufficiently small and $0 < \de'' < \de'$, there is an $r > 0$ such that the following holds. Let $a \in B(0,r)$ and assume that $|\xi_k(a)-\mu_k(a)| \leq \de'$, for all $k \leq n$ and $|\xi_n(a)-\mu_n(a)| \geq \de''$. Then
\[
\biggl| \frac{\xi_n'(a)}{(R_a^n)'(\mu_0(a))x'(a)} - 1 \biggr| \leq \vep.
\]
\end{Lem}
\begin{proof}
First we note that by Lemma \ref{smallep} we have
\[
\xi_n(a) = x(a) (R_a^n)'(\mu_0(a)) + \mu_n(a) + E_1(a),
\]
where $|E_1(a)| \leq |\xi_n(a) - \mu_n(a)|/1000$ independently of $n$ and $a$. Put $R_a'(\mu_j(a)) = \la_{a,j}$. Differentiating with respect to $a$ we get
\begin{equation} \label{xip}
\xi_n'(a) = \prod_{j=0}^{n-1} \la_{a,j} \biggr( x'(a) + x(a) \sum_{j=0}^{n-1} \frac{\la_{a,j}'}{\la_{a,j}} + \frac{\mu_n'(a) + E_1'(a)}{\prod_{j=0}^{n-1} \la_{a,j}}  \biggl).
\end{equation}

We claim that only the $x'(a)$ is dominant in (\ref{xip}) if $n$ is large so that $\de'' \leq |\xi_n(a)-\mu_n(a)| \leq \de'$. This means that, again by Lemma \ref{smallep},
\[
\de'' \lesssim |x(a)| \prod_{j=0}^{n-1} |\la_{a,j}| \lesssim \de' < 1.
\]
Since $\prod_{j=0}^{n-1} |\la_{a,j}| \geq \la^n$, for some $\la > 1$, taking logarithms and rearranging we get
\begin{equation}
n \sim - \frac{\log |x(a)|}{\frac{1}{n} \sum_{j=0}^{n-1} \log |\la_{a,j}|} \text{    or     } -\log |x(a)| \sim \sum_{j=0}^{n-1} \log |\la_{a,j}|, \label{nsize}
\end{equation}
if $|\log \de''|  \ll |\log |x(a)||$, which is true if the perturbation $r > 0$ is chosen sufficiently small compared to $\de''$.
Since $|\la_{a,j}| \geq \la > 1$, this means that
\[
|x(a)| \sum_{j=0}^{n-1} \frac{|\la_{a,j}'|}{|\la_{a,j}|} \leq |x(a)| n C \leq -C |x(a)| \log|x(a)|.
\]
Finally $-|x(a)| \log|x(a)| / |x'(a)| \raw 0$ as $a \raw 0$.

The last two terms in (\ref{xip}) tend to zero as $n \raw \infty$, since $\mu_n'(a)$ and $E_1'(a)$ are bounded. We have proved that
\[
\xi_n'(a) \sim x'(a) \prod_{j=0}^{n-1} \la_{a,j},
\]
if $|\xi_n(a)-\mu_n(a)| \leq \de'$ and $n \geq N$ for some $N$. Choose the perturbation $r$ sufficiently small so that this $N$ is at most the number $n$ determined by (\ref{nsize}). Since $\la_{a,j} = R_a'(\mu_j(a))$, the proof is finished.
\end{proof}

The following lemma, which will be needed in the subsequent lemma, is variant of Lemma 15.3 in \cite{WR} (see also \cite{MA}, Lemma 2.1).
\begin{Lem} \label{prod-dist}
Let $u_n \in \C$ be complex numbers for $1 \leq n \leq N$. Then
\begin{equation}
\biggl| \prod_{n=1}^N (1+u_n)-1 \biggr| \leq
\exp \biggl(\sum_{n=1}^N |u_n| \biggr) - 1.
\label{ineq-1}
\end{equation}
\end{Lem}

\begin{Lem}[Distortion lemma] \label{distortion1}
For every $\vep > 0$, there are arbitrarily small constants $\de' > 0$ and $r > 0$ such that the following holds. Let $a,b \in B(0,r)$ and suppose that $|\xi_k(t) - \mu_k(t)| \leq \de'$, for $t=a,b$ and all $k \leq n$. Then
\[
\biggl| \frac{(R^n)'(v(a),a)}{(R^n)'(v(b),b)} - 1 \biggr| < \vep.
\]
\end{Lem}
The same statement holds if one replaces $v(s)=\xi_0(s)$, $s=a,b$, by $\mu_0(t)$, $t=a,b$.
\begin{proof}
The proof goes in two steps. Let us first show that
\begin{equation}
\biggl| \frac{(R_t^n)'(\mu_0(t))}{(R_t^n)'(\xi_0(t))} - 1 \biggr| \leq \vep_1,
\label{close-to-1}
\end{equation}
where $\vep_1=\vep(\de')$ is very close to $0$. We have
\begin{align}
\sum_{j=0}^{n-1} \biggl| \frac{R_t'(\mu_j(t)) - R_t'(\xi_j(t))}{R_t'(\xi_j(t))} \biggr| &\leq C \sum_{j=0}^{n-1} |R_t'(\mu_j(t)) - R_t'(\xi_j(t))| \nonumber \\
&\leq C \sum_{j=0}^{n-1} |\mu_j(t)-\xi_j(t)| \nonumber \\
&\leq C \sum_{j=0}^{n-1} \la^{j-n}|\mu_n(t)-\xi_n(t)| \leq C(\de'), \nonumber
\end{align}
where we used Lemma \ref{expinside}. By Lemma \ref{prod-dist}, (\ref{close-to-1}) holds of $\de'$ is small enough.
Secondly, we show that
\[
\biggl| \frac{(R_t^n)'(\mu_0(t))}{(R_s^n)'(\mu_0(s))} - 1\biggr|\leq \vep_2,
\]
where $\vep_2=\vep_2(\de')$ is close to $0$. Since $\la_{t,j}$ are all analytic in $t$ we have
$\la_{t,j}=\la_{0,j}(1+c_j t^l + \ldots)$. Moreover, since $n \leq -C \log|x(t)| = -C \log |t|$
\[
\frac{(R_t^n)'(\mu_0(t))}{(R_s^n)'(\mu_0(s))}  = \prod \frac{\la_{t,j}}{\la_{s,j}} = \prod_{j=0}^{n-1}
\frac{\la_{0,j}(1+c_j t^l + \ldots)}{\la_{0,j}(1+c_j s^l + \ldots)} =
\frac{1+cn t^l + \ldots}{1+cn s^l + \ldots}.
\]
Both the last numerator and denominator in the above equation can be
estimated by $1 + c' (\log |t|) |t|^l$ and $1 + c' (\log |s|) |s|^l$. The lemma follows easily.
\end{proof}

Combining Lemma \ref{xinprim} and Lemma \ref{distortion1} we immediately get
\begin{Cor} \label{vava}
Let $a \in B(0,r)$, $0 < \de'' < \de'$ and assume that $|\xi_k(a)-\mu_k(a)| \leq \de'$, for all $k \leq n$ and $|\xi_n(a) - \mu_n(a)| \geq \de''$. Then
\[
\xi_n'(a) \sim (R_a^n)'(v(a)) x'(a).
\]
\end{Cor}
In the following lemma assume always that the disk $D_0$ satisfies $D_0=B(a_0,k_0|a_0|)$, where $k_0$ is such that $a,b \in D_0$ implies $|x(a)/x(b) - 1| \leq 1/10$, (see (\ref{xK1})).
\begin{Lem}\label{lemma36}
For $a,b \in D_0$ and $0 < \de'' < \de'$, assume that $|\xi_k(a)-\mu_k(a)| \leq \de'$ for all $k \leq n$ and 
and $|\xi_j(a)-\mu_n(a)| \geq \de''$. Then
\[
 |\xi_n(a)-\xi_n(b)| \sim \prod_{j=0}^{n-1} |\la_{0,j}||x(a)-x(b)|.
\]
\end{Lem}
\begin{proof}
The condition on $D_0$ and $|\xi_n(a)-\mu_n(a)| \leq \de'$ implies by equation (\ref{xK1}), Lemma \ref{distortion1} and Lemma \ref{smallep} that $|\xi_n(b)-\mu_n(b)| \leq 2\de'$ and $\xi_k(b) \in \NN$, for all $k \leq n$. It also implies that $n \leq -C \log |x(a)|$ (cf. (\ref{nsize})). We have
\begin{align}
|\xi_n(a)-\xi_n(b)| &\geq |\xi_n(a)-\mu_n(a)-(\xi_n(b)-\mu_n(b))|-|\mu_n(a)-\mu_n(b)|
\nonumber \\
&\sim \bigl|x(a) \prod_{j=0}^{n-1} \la_{a,j}-x(b) \prod_{j=0}^{n-1} \la_{b,j} \bigr| - |\mu_n(a)-\mu_n(b)| \nonumber \\
&\sim \prod_j |\la_{0,j}| \bigl| x(a) \prod_j (1+c_j a^{l_j}) - x(b) \prod_j (1+c_j b^{l_j}) \bigr| -|\mu_n(a)-\mu_n(b)| \nonumber \\
&\sim \prod_j |\la_{0,j}| | x(a)(1+n c a^l) - x(b)(1+n c b^l)| -
|\mu_n(a)-\mu_n(b)| \nonumber \\
&\sim \prod_j |\la_{0,j}||x(a)-x(b)| - |\mu_n(a)-\mu_n(b)|. \label{firstterm}
\end{align}

We want the first term in (\ref{firstterm}) to be dominant over the second term. The condition $\de'' \leq |\xi_n(a)-\mu_n(a)|$ implies $-\log |x(a)|  \lesssim \sum_j \log |\la_{a,j}|$ (cf. (\ref{nsize})). By Lemma \ref{distortion1}, $\sum_j \log |\la_{a,j}| \sim \sum_j \log |\la_{0,j}|$. The condition on $D_0$ means that $1/C \leq |x'(a)||a-b|/|x(a)-x(b)| \leq C$, for some $C > 1$, ($C \raw 1$ as $k_0 \raw 0$). Therefore, we can estimate the first term in (\ref{firstterm}) by
\begin{align}
\prod_{j=0}^{n-1} |\la_{0,j}||x'(a)||a-b| &\sim \exp ( (\sum_{j=0}^{n-1} \log |\la_{0,j}|) + \frac{k-1}{k} \log |x(a)|) |a-b| \nonumber \\
&\gtrsim \exp ( (\sum_{j=0}^{n-1} \log |\la_{0,j}|)/k ) |a-b| \nonumber \\
&\gtrsim e^{\ga' n} |a-b|, \nonumber
\end{align}
where $k \ga' = \log |\la|$, and $\la$ is as in Lemma \ref{expinside}.
Since $|\mu_n'(a)|$ is bounded, this means that
\[
|\xi_n(a)-\xi_n(b)| \gtrsim \prod_j |\la_{0,j}||x(a)-x(b)|.
\]
The other inequality follows in precisely the same way.
\end{proof}


\begin{proof}[Proof of Proposition \ref{da/dz}]


First, we will use Lemma \ref{xinprim} and prove by induction, that
\begin{equation}
|\xi_{N+k}'(a)| \geq e^{\ga'(N+k)}, \label{initexp}
\end{equation}
where $\ga' = \ga/k - \vep_0$, for some small $\vep_0 > 0$, $\ga = \log \la$, $\la$ is as in Lemma \ref{expinside}, and $k$ is as in (\ref{xK1}). Let $N$ be the smallest integer such that $|\xi_N(a)-\mu_N(a)| \geq \de''$. Equation (\ref{initexp}) is fulfilled for $k=0$, if $N$ is sufficiently large, hence if $B(0,r)$ is sufficiently small. Indeed, by Corollary \ref{vava} and Lemma \ref{expinside}, we have
\begin{equation}
|\xi_N'(a)| \geq (1/2)|(R_a^N)'(v(a))||x'(a)| \geq (1/4) e^{\ga
  N}|x'(a)| \geq e^{\ga' N}.
\label{gaprim}
\end{equation}
Hence the initial condition ($k=0$) of the induction is satisfied.

So, assume that
\[
|\xi_{N+j}'(a)| \geq e^{\ga'(N+j)},  \text{ for all $j\leq k$}.
\]
We want to prove that
\[
|\xi_{N+j}'(a)| \geq e^{\ga'(N+j)},  \text{ for all $j\leq k+1$}.
\]
First note that the first assumption on $a$ gives
\begin{equation}
|R'(\xi_j(a),a)| \geq C_1^{-1},
\label{bbaa}
\end{equation}
where $C_1 = C_1(\de)$.

By the Chain Rule we have the recursions (remember the notation
$\xi_{n}(a) = R^n(v(a),a)$)
\begin{align}
\frac{\partial R^{n+1}(v(a),a)}{\partial z} &= \frac{\partial
  R(\xi_n(a),a)}{\partial z} \frac{\partial R^n(v(a),a)}{\partial z}
  \label{dz}, \\
\frac{\partial R^{n+1}(v(a),a)}{\partial a} &= \frac{\partial
  R(\xi_n(a),a)}{\partial z} \frac{\partial R^n(v(a),a)}{\partial a} +
\frac{\partial R(\xi_n(a),a)}{\partial a}. \label{da}
\end{align}
Now, the recursion formulas (\ref{dz}) and (\ref{da}), together
with (\ref{bbaa}) and (\ref{gaprim}), gives
\begin{align}
|\xi_{N+k+1}'(a)| &\geq
 |R_a'(\xi_{N+k}(a))||\xi_{N+k}'(a)|\biggl(1-\frac{|\partial_a
   R_a(\xi_{N+k}(a))|}
{|R_a'(\xi_{N+k}(a))||\xi_{N+k}'(a)|} \biggr) \nonumber \\
&\geq
|(R_a^{k+1})'(\xi_N(a))||\xi_N'(a)| \prod_{j=0}^k \biggl(1-C_1 \frac{|\partial_a
   R_a(\xi_{N+j}(a))|}{|\xi_{N+j}'(a)|} \biggr)
\nonumber \\
&\geq
e^{\ga (k+1)} e^{\ga' N} \prod_{j=0}^k (1-BC_1 e^{-\ga'(N+j)} ) \nonumber \\
&\geq e^{(\ga-\ga')(k+1)} e^{\ga' (N + k + 1)}
\prod_{j=0}^k (1-B' e^{-\ga'(N+j)} ) \geq e^{\ga' (N+k+1)},  \nonumber
\end{align}
if $N$ is large enough, (here $B' = BC_1$). The sum
\[
\sum_{j=0}^{\infty} B'e^{-\ga'(N+j)} < \infty,
\]
and can be made arbitrarily small if $N$ is large enough.

By the definition of $Q_n(a)$, we have
\[
Q_{N+n}(a) = Q_N(a) \prod_{j=0}^{n} \biggl(1+\frac{\partial_a
   R_a(\xi_{N+j}(a))}{R_a'(\xi_{N+j}(a)) \xi_{N+j}'(a)} \biggr).
\]
So,
\begin{equation}
|Q_{N+n}(a)-Q_N(a)| \leq  |Q_N(a)|/1000, \nonumber
\end{equation}
if $N$ is sufficiently large.
\end{proof}


\begin{Rem} \label{QN}
The number $Q_N(a)$, for general $a \in D_0$, can be estimated
by Corollary \ref{vava} in the following way;
\begin{equation}
\biggl| \frac{Q_N(a)}{x'(a)} - 1 \biggr| \leq \vep_1,
\nonumber
\end{equation}
where $\vep_1 > 0$ can be made arbitrarily small if $r > 0$ is small enough.
If we want good argument distortion of a small disk $D_0 = B(a_0,k_0|a_0|)$,
i.e. the quotient $Q_N(a)/Q_N(a_0)$ is very close to $1$ for all $a
 \in D_0$, then we must have
\begin{equation}
\biggl| \frac{x'(a)}{x'(b)} -1 \biggr| \leq \vep_2,
\label{qna}
\end{equation}
for all $a \in D_0$, for some $\vep_2 > 0$ small enough.
Equation (\ref{qna}) gives an estimate of the number $k_0$; it follows from
(\ref{xK1}) that it is enough to have
\begin{equation}
|k_0^{k-1}|\leq \vep_3,
\label{k0}
\end{equation}
for some suitable $\vep_3 = \vep_3(\vep_2)$.
\end{Rem}

\begin{Cor} \label{QN-cor}
If $D_0 = B(a_0,k_0|a_0|)$, where $k_0$ satisfies (\ref{k0}), then
\begin{equation}
\biggl| \frac{Q_{N}(a)}{Q_{N}(a_0)} - 1 \biggr| \leq 1/500,
\label{QNa0}
\end{equation}
for all $a \in D_0$, (if the $\vep_j$:s are chosen suitable).
\end{Cor}

\begin{Rem}
We also see that
\begin{equation}
\frac{1}{2} |x'(a)| \leq |Q_N(a)| =
\frac{|\xi_N'(a)|}{|(R^N)'(v(a),a)|}
\leq 2 |x'(a)|,
\nonumber
\end{equation}
if $a \in D_0$. In particular, if $x'(0) \neq 0$ then $Q_N(a) = K_1 + \OO(a)$ in a neighborhood of $a=0$,
so the equation (\ref{QNa0}) is valid in a whole disk $B(0,r)$
instead of only a small $D_0$ "far" away from $0$ (meaning $r_0 \ll |a_0|$).
The parameter directions for which $x'(0) \neq 0$ are usually called {\em non-degenerate}, (see also \cite{MA}).
\end{Rem}


In the next lemma we show that the disk $D_0=D(a_0,r_0)$ will grow to size at least $S$ before it leaves $\NN$. Note that $S > 0$ depends only on $\de'$.

\begin{Lem} \label{initdist}
If $r > 0$ is sufficiently small then there are numbers $k_0$, $0 < k_0 < 1$ and $S=S(\de')$, such that if we chose $r_0 = k_0|a_0|$ then the following holds for any dyadic disk $D_0=B(a_0,r_0) \subset B(0,r)$:
There is an $n$ such that the set $\xi_{n}(D_0) \subset P_{\de'}$ and has diameter at least $S$. Moreover, we have low argument distortion, i.e.
\begin{equation}
\biggl| \frac{\xi_n'(a)}{\xi_n(b)} - 1 \biggr| \leq 1/100,
\label{xinn}
\end{equation}
for all $a,b \in D_0$.
\end{Lem}
\begin{proof}
Choose the maximal $n$ such that $\xi_k(a_0) \in \NN$ for all $k \leq n$ and \[
\de'/(2M_0) \leq |\xi_n(a_0)-\mu_n(a_0)| \leq \de'/2,
\]
where $M_0$ is the supremum of $|R_a'(z)|$ over all $a \in B(0,r)$ and $z \in \hat{\C}$.

Putting $a=a_0$ and $b=a$ in Lemma \ref{lemma36}, we see that $|\xi_n(a)-\xi_n(a_0)|$ expands almost linearly for all $a \in D_0=B(a_0,k_0|a_0|)$. Choose $k_0 > 0$ maximal, subject to the condition preceding Lemma \ref{lemma36}, such that the diameter of $\xi_n(D_0)$ is at most $\de'/(10M_0)$. Then every $b \in D_0$ has that
\[
\de'' \leq |\xi_n(b)-\mu_n(b)| \leq \de'.
\]

Now, impose another condition on $k_0$ such that (\ref{xinn}) is fulfilled
for all $a,b \in D_0 = D(a_0,k_0 |a_0|)$ (this possibly diminishes $k_0$).
Combining Lemma \ref{distortion1}, Proposition \ref{da/dz} and Corollary \ref{QN-cor}, we have that
(\ref{xinn}) is fulfilled if $|k_0^{k-1}| \leq \vep$, for some sufficiently small $\vep > 0$.

By (\ref{xinn}) and Lemma \ref{xinprim}
the diameter $D$ of the set $\xi_n(D_0)$ can be estimated by
\begin{align}
D &\geq |\xi_n'(a_0)||a_0k_0| \geq (1/2) \prod_j |\la_{a_0,j}||x'(a_0)||a_0 k_0|
\geq C \prod_j |\la_{a_0,j}||a_0|^{k}k_0.
\nonumber
\end{align}

By Lemma \ref{distortion1}, it follows that
\begin{equation}
\biggl| \frac{\prod_j \la_{a,j}}{\prod_j \la_{b,j}} - 1\biggr| \leq \vep,
\label{1st-dist}
\end{equation}
for all $a,b \in B(0,r)$ if $r$ is small enough.
By Lemma \ref{smallep} and (\ref{1st-dist}),
\begin{align}
\de'/(2M_0) &\leq |\xi_n(a_0)-\mu_n(a_0)| \leq C \prod_j |\la_{0,j}||x(a_0)|
=C \prod_j |\la_{0,j}||a_0|^k.
\end{align}
Thus,
\begin{equation}
\frac{D}{\de'} \geq C \frac{\prod_j |\la_{a,j}|}{\prod_j |\la_{0,j}|}k_0 \geq
C k_0.\nonumber
\end{equation}
So the diameter $D$ of the set $\xi_n(D_0)$ is greater than $S = S(\de')$.
Also, by (\ref{xinn}), we have bounded argument distortion for all $a,b \in D_0$.
\end{proof}

The next lemma provides strong bounded distortion as long as two points in $\NN$ stays at bounded distance from each other under iteration and if we admit a finite number $\ti{N}$ of iterates and not meeting $U_{\de/10}$.
\begin{Lem}[Extended Distortion Lemma]\label{distortion2}
Let $\ti{N} \in \N$ and $\vep > 0$. There exists an $r > 0$ and $S_1 > 0$ such that the following holds. Let $a,b \in B(0,r)$ and assume that $z,w \in \NN$ are such that $R^k(z,a), R^k(w,b) \notin U_{\de /10}$ and $|R^k(z,a)-R^k(w,b)| \leq S_1$ for all $k=0, \ldots,n$, where $n \leq \ti{N}$.
Then
\[
\biggl| \frac{(R^n)'(z,a)}{(R^n)'(w,b)} -1 \biggr| < \vep.
\]
\end{Lem}
\begin{proof}
Put $R_a^j(z)=z_j$ and $R_b^j(w)=w_j$. Since the parameter dependence can be made arbitrarily small under $\ti{N}$ iterations by choosing $r > 0$ sufficiently small, using Lemma \ref{prod-dist} we estimate the sum
\[
\sum_{j=0}^n \frac{|R_a'(z_j)-R_b'(w_j)|}{|R_b'(w_j)|} \leq C(\de) \sum_{j=0}^n |z_j-w_j| \leq C(\de) \ti{N} S_1.
\]
Letting $S_1 \leq \vep' / (C(\de) \ti{N})$, for $\vep' = \log (1+\vep)$, we get the desired result.
\end{proof}

We also need a global distortion lemma, which is valid if we go over the scale $S_1$, but only admit a finite
number $\ti{N}$ of iterates and not meeting $U_{\de/10}$. Here we relax the condition $|R^k(z,a)-R^k(w,b)| \leq S_1$ for $k \leq \ti{N}$, for some $\ti{N} \in \N$.
In this case the distortion is not necessarily so low anymore, but still bounded. The following is immediate.
\begin{Lem}[Global Distortion Lemma] \label{distortion3}
Let $\ti{N}$ and $r > 0$ be as in the above lemma and $a,b \in B(0,r)$. Assume that $z,w$ are such that $R^k(z,a), R^k(w,b) \notin U_{\de/10}$ for $k \leq \ti{N}$. Then
\[
\biggl| \frac{(R^n)'(z,a)}{(R^n)'(w,b)}\biggl| < \ti{C}(\ti{N}).
\]
\end{Lem}


\section{The free period and the degree of $\xi_n$}
The main object of this section is to show that once the set
$\xi_n(D_0)$ has reached diameter $S=S(\de') > 0$ (which follows by Lemma \ref{initdist}) then
$\xi_{n+m}(D_0)$ will intersect $U_{\de/10}$ within a finite number of
iterates, i.e. $m \leq \ti{N}$ for some $\ti{N}$ only depending on $\de'$. These last $m$ iterates are referred to as the free period.


\begin{Lem} \label{finret}
For every $d > 0$ there is an $r > 0$ such that the following holds. If $z \in J(R_{a}) \cap P_{\de'}$, for some $a \in B(0,r)$, then there is a disk $D$ of diameter $d$ containing $z$ and a positive integer $\ti{N}$, which only depends on $d$, such that
\[
\sup_{a \in B(0,r)} \inf \{ m \in \N : R_{a}^m(D) \cap U_{\de/10} \neq \emptyset \} \leq \ti{N}.
\]
\end{Lem}

\begin{proof}
First cover the closure of $P_{\de'}$ with a finite collection of open disks $D_j$ of diameter $d$. Assume further that every disk $D_j \subset F(R_0)$ is compactly contained in $F(R_0)$.
Let $D_j$, $j \in I$, be the subcollection of disks such that $D_j \cap J(R_0) \neq \emptyset$.
Since $R_0^n$ is not normal on the Julia set, we get that for every $D_j$, $j \in I$, there is a smallest number $n=n(j)$ such that $R_0^n(D_j) \cap U_{\de/10} \neq \emptyset$. Since $R_a^n(D_j)$ moves continuously in $a$, there is an $r > 0$ such that the same statement holds for $R_a$ instead of $R_0$, if $a \in B(0,r)$ and such that those $D_j \Subset F(R_0)$ still have $D_j \Subset F(R_a)$ (recall that $F(R_0)$ corresponds only to super-attracting cycles).
\end{proof}

If $d=S/2$, note that $\ti{N}$ depends only on $\de'$, since $S=S(\de')$.
In the following assume always that $D_0 \subset B(0,r)$ is a disk $D_0=(a_0,r_0)$ which has that $r_0=|a_0|k_0$, where $k_0$ is fixed and satisfies (\ref{k0}).

\begin{Lem} \label{return}
Assume that $\xi_n(D_0)$ has diameter at least $S$ and $\xi_n(D_0) \subset P_{\de'}$, where $n$ is chosen as in Lemma \ref{initdist}. Then there is a maximum number $\ti{N} < \infty$ such that either
\[
\sup_{D_0 \subset B(0,r)} \inf \{ m \in \N : \xi_{n+m}(D_0) \cap U_{\de/10} \neq \emptyset \} \leq \ti{N},
\]
or there is a topological disk $D_1 \subset D_0$ such that $\xi_n(D_1)$ has diameter at least $S/4$ and such that for all but countably many $a \in D_1$, $R_a$ is not a Misiurewicz map.
\end{Lem}
\begin{proof}
Let the diameter of $D=\xi_n(D_0)$ be $d \geq S$. Consider two cases. Either there is a ball $D' \subset D$, with the same center as $D$ and radius $S/4$, which is contained in $F(R_{a})$ for all $a \in D_1=\xi_n^{-1}(D')$ or not. If the former case occurs then all but countably many parameters in $D_1$ corresponds to maps which are not Misiurewicz maps, (the only exception are those parameters which are centers of hyperbolic components). If $D' \cap J(R_{a}) \neq \emptyset$ for some $a \in D_1$, then it means that any ball $D'' \subset D$ with radius $S/2$ containing a point $z \in D' \cap J(R_{a})$ has that $D'' \subset D$. Now use Lemma \ref{finret} on some ball $D''$, i.e. $d=S/2$.

We get an upper bound of the return time of $R_{a}^m(D)$. To switch from $R_{a}^m(D)$ to $\xi_{n+m}(D_0)$ we note that, precisely as in Lemma \ref{finret}, the return time is locally constant for small perturbations of parameters.
\end{proof}

Thus, we have proved that if a disc $D_0$ grows to size $S$ then a
return will occur into $U_{\de/10}$ after at most $\ti{N}$ iterates, provided the perturbation is sufficiently small. Since $\xi_n$ is injective until it reaches size $S$, and since there is a finite number of iterates after this time until a return occurs into $U_{\de/10}$, it is not hard to see that the degree of $\xi_n(a)$ is bounded for all $a \in D_0$, uniformly for all $D_0 \subset B(0,r)$. This is the content of the next lemma.
\begin{Lem}
Assume that the diameter of $\xi_n(D_0)$ is at least $S$ for the smallest possible $n$, and suppose that $\xi_n(D_1) \cap J(R_a)$ for some $a \in D_1$ where $D_1$ is as in Lemma \ref{return}, (i.e. the second case in this lemma does not occur). Let $m$ be minimal so that $\xi_{n+m}(D_0) \cap U_{\de/10} \neq \emptyset$. Then the degree of $\xi_{n+m}$ on $D_0$ is bounded by some $M < \infty$, regardless of $n$.
\end{Lem}

\begin{proof}
Consider the almost round disk $D=\xi_n(D_0)$ centered at $z \in P_{\de'} \cap J(R_a)$, ($D$ is almost round, by Lemma \ref{initdist}). The degree of $R_0^m$ on $D$ is bounded by $d^{\ti{N}}$. Moreover, for any point $z \in D$, by definition we have that $R_a^j(z)$ does not intersect $U_{\de/10}$ for $j=0,\ldots, m-1$. This means that $|(R_a^m)'(z)| \geq C $ for some $C=C(\de,\ti{N})$ and there are no critical points of $R_a^m$ inside $D$. It follows that any two points $z_1,z_2 \in D$ mapped onto the same point must be separated by some fixed constant $c$. Now, to switch from $R_a^m(D)$ to $\xi_{n+m}(D_0)$ we note that each pair of points $z_1,z_2 \in D$ are images under $\xi_n(a)$ for some $a \in D_0$, i.e. $z_1 = \xi_n(a_1), z_2 = \xi_n(a_2)$. The parameter dependence under the coming $m \leq \ti{N}$ iterates can be made arbitrarily small if the perturbation is sufficiently small (i.e if $r > 0$ is sufficiently small). Hence there is a slightly smaller constant $c/2$ such that if $z_1=\xi_n(a_1), z_2 = \xi_n(a_2)$ then $\xi_{n+m}(a_1)=\xi_{n+m}(a_2)$ only if $|z_1-z_2| \geq c/2$ or $z_1=z_2$. The lemma follows.
\end{proof}

\section{Conclusion and proof of Theorem \ref{discthm}}
We want to show that a specific fraction $f > 0$ of parameters in
the disc $B(0,r)$ are not $(\de,k)$-Misiurewicz for any $k \geq 0$, by showing that either there is a specific fraction that corresponds to those maps, for which a critical point returns into $U_{\de}$ after $k$ iterations or simply that these parameters correspond to a hyperbolic component (this corresponds to the second case in Lemma \ref{return}). Choosing $r$ sufficiently small and $n$ minimal in Lemma \ref{initdist}, so that $\xi_n(D_0)$ has diameter at least $S$ for $n > k$, in the second case we have immediately that the fraction of parameters in $D_0$, which do not correspond $(\de,k)$-Misiurewicz maps, is bounded by some number $f > 0$ (not depending on $\de$). Let us therefore focus on the first case.

Let us assume that $n=N_1+m$ is smallest integer such that $\xi_{n}(D_0) \cap U_{\de/10} \neq \emptyset$, where $N_1$ is the minimal number of iterates for which that diameter of $\xi_{N_1}(D_0) \subset \NN$ is at least $S$. Moreover, for any $k$ we can choose the perturbation $r$ sufficiently small so that $N_1 > k$. Hence any rational map $R_a$, for $a \in B(0,r)$, which has a critical point $c(a)$ not in a super-attracting cycle and such that $c(a)$ returns into a slightly smaller $U_{(9/10)\de} \subset U_{\de}$ is not $(\de,k)$-Misiurewicz.

Choosing $N_1$ sufficiently large (i.e. $r$ sufficiently small) we can ensure that for all $a \in D_0$,
\begin{equation}
|(R_a^{N_1+k})'(v(a))| \geq e^{\ga (N_1 + k)},
\label{expoutside*}
\end{equation}
for all $k \leq m \leq \ti{N}$, for some $\ga \geq \uli{\ga}$, ($\uli{\ga} = \log \la/2$, where $\la$ is as in Lemma \ref{expinside}).

\begin{Rem}
Equation (\ref{expoutside*}) can also be shown to hold at the return time $k=m$ by using a variant of the  Outside Expansion Lemma (Lemma 3.7 in \cite{MA}), which says that $|(R_a^m)'(z)| \geq C \la^m$, for some $\la > 1$ as long as $R_a^k(z) \notin U_{\de/10}$, for all $k \leq m-1$ and $R_a^m(z) \in U_{\de/10}$.
\end{Rem}

Let $N \leq N_1$ in Proposition \ref{da/dz}. We get
\[
\biggl| \frac{Q_N(a) (R^n)'(v(a),a)}{\xi_n'(a)} - 1 \biggr| \leq 1/1000,
\]
for all $a \in D_0$. 
Now the Distortion Lemmas \ref{distortion1} and \ref{distortion2} implies that as long as $|\xi_k(a) - \xi_k(b)| \leq S_1$, for all $k \leq n$, we have
\begin{equation}
\biggl|\frac{(R_a^n)'(v(a))}{(R_b^n)'(v(b))} - 1 \biggr| \leq \vep,
\nonumber
\end{equation}
for some $\vep > 0$. 
Hence by Corollary \ref{QN-cor}, choosing $\vep > 0$ suitable, the following holds:
\begin{equation}
\biggl| \frac{\xi_{n}'(a)}{\xi_{n}'(b)} - 1 \biggr| \leq 1/100, \label{angledist}
\end{equation}
for all $a,b \in D_0$ as long as $|\xi_n(a)-\xi_n(b)| \leq S_1$.
This is the geometry control we need.

By the strong distortion estimate up to scale $S_1$, (Lemma \ref{distortion1} and Lemma \ref{distortion2}) together with weaker but global Distortion Lemma \ref{distortion3} it follows that there is a constant $C_1$ such that if $a,b \in D_0$ and if $z,w \in D=\xi_{N_1}(D_0)$, then
\[
\biggl| \frac{(R_a^k)'(z)}{(R_b^k)'(w)} \biggr| \leq \ti{C},
\]
if $k \leq m$, where $m$ is the first return time into $U_{\de/10}$. By (\ref{expoutside*}) and Proposition \ref{da/dz} we have that
\[
\biggl| \frac{\xi_{N_1+k}'(a)}{\xi_{N_1+k}'(b)} \biggr| \leq C,
\]
for $k \leq m$, where $C=C(\ti{N})$, (recall that $m$ is bounded by $\ti{N}$). Since $\xi_{n}$ is at most $M$-to-$1$ on $D_0$, this means that those parameters mapped into $U_{(9/10)\de}$ correspond to a definite fraction of the part of $\xi_{n}(D_0)$ that cover $U_{(9/10)\de}$. If $\mu$ is the Lebesgue measure, we get
\[
\frac{\mu (\{ a \in D_0 : \xi_{n}(a) \in U_{(9/10)\de}  \}) }{\mu(D_0)} \geq C_0 \frac{E}{F}
\]
where $E$ is the area of $\xi_{n}(D_0) \cap U_{(9/10)\de}$, $F$ is the area of $\xi_n(D_0)$ and $C_0$ is a constant depending on the degree of $\xi_{n}$ (which is bounded) and $\de$.

Since $\xi_{n}(D_0)$ intersects $U_{\de/10}$, the strong distortion estimate (\ref{angledist}) inside $D_0$
implies that the area $E$ is at least the area of $U_{(9/10)\de} \cap C$ where $C$ is a circle with radius $r_1 \de$, where $r_1$ only depends on $\de$. In particular, $\mu(E) \geq C_1$, where $C_1=C_1(\de)$. Since the area of $\xi_{n}(D_0)$ is bounded by the area of the Riemann sphere, which is bounded by some $C_2 < \infty$, we get
\[
\frac{\mu (\{ a \in D_0 : \xi_{n}(a) \in U_{\de}  \})}{\mu(D_0)} \geq C_0 \frac{C_1}{C_2} \equiv f > 0.
\]
Since this estimate holds for every small disk $D_0 \subset B(0,r)$,
we conclude that the set of parameters $a \in B(0,r)$ of
$(\de,k)$-Misiurewicz maps has Lebesgue density at most $1-f < 1$, at
$a=0$. The proof of Theorem \ref{discthm} is finished.


\bibliographystyle{plain}
\bibliography{ref}

\comm{
\bibliographystyle{plain}
\bibliography{bib}

}

\end{document}